# Bayesian Inference For Exponential Distribution Based On Upper Record Range

**P. Nasiri, S. Hosseini, M. Yarmohammadi
and F. Hatami**

**Abstract** This paper deals with Bayesian estimations of scale parameter of the exponential distribution based on upper record range $(R_n)$. This has been done in two steps; point and interval. In the first step the quadratic, squared error and absolute error, loss functions have been considered to obtain Bayesian-point estimations. Also in the next step the shortest Bayes interval (Hight Posterior Density interval) and Bayes interval with equal tails based on upper record range have been found. Therefore, the Homotopy Perturbation Method (HPM) has been applied to obtain the limits of Hight Posterior Density intervals. Moreover, efforts have been made to meet the admissibility conditions for linear estimators based on upper record range of the form $mR_n + d$ by obtained Bayesian point estimations. So regarding the consideration of loss functions, the prior distribution between the conjunction family has been chosen to be able to produce the linear estimations from upper record range statistics. Finally, some numerical examples and simulations have been presented.

P. Nasiri
Department of Statistics, University of Payame Noor, 19395-4697 Tehran, I. R. of IRAN
E-mail: pnasiri45@yahoo.com

S. Hosseini
Department of Statistics, Technical Institute of Koye, Koye, Kurdistan Region,IRAQ
E-mail: S.Hosseini.stat@gmail.com

M. Yarmohammadi
Department of Statistics, University of Payame Noor, 19395-4697 Tehran, I. R. of IRAN
E-mail: masyar@pnu.ac.ir

F. Hatami
Department of Mathematics, Iran University of Science and Technology, Tehran, I. R. of IRAN
E-mail: Farhad.hatami67@gmail.com

# 1 Introduction

Let $X_1, X_2, X_3, \ldots$ be a sequence of independent and identically distributed (iid) random variables with cumulative distribution function (CDF) F(x) and probability density function (pdf) f(x). For $n \geq 1$, define:

$$T(1) = 1, T(n+1) = min\{j : X_j \geq X_T(n)\},$$

the sequence $\{X_{T(n)}\}_{n=1}^{\infty}$ is known as upper record values and the sequence $\{T(n)\}_{n=1}^{\infty}$ is known as record times sequence (Arnold et al. (1998)). Chandler (1952), was the first who defined and framed the concepts of record values, record times and related statistics theoretically. Interested readers may refer to Arnold (1998), Nagaraja (1988) and Nevzorov (1946), for further information on fundamental concepts in records. There are scholars who have provide some pure and inferential studies based on record values and record times, among whom Ahsanullah (1990), Balakrishnan, Ahsanullah, Chan (1995), Feuerverger and Hall (1998) are the most famous ones.

Suppose X is an exponential random variable with the parameter $\delta$, which its probability density function (PDF) and Cumulative distribution function (CDF) are respectively as followed:

$$f(x;\delta) = \frac{1}{\delta} e^{-(\frac{x}{\delta})}, x \geq 0, \delta > 0, \tag{1}$$

and

$$F(x;\delta) = 1 - e^{-(\frac{x}{\delta})}, x \geq 0, \delta > 0. \tag{2}$$

The exponential distribution is enormously useful in different contexts of science and Technology. The exponential distribution occurs naturally when describing the lengths of the inter-arrival times in a homogeneous Poisson process. Exponential variables can also be used to model situations where certain events occur with a constant probability per unit length, such as the distance between mutations on a DNA strand, or between road-kills on a given road. In queuing theory, the service times of agents in a system (e.g. how long it takes for a bank teller etc. to serve a customer) are often modeled as exponentially distributed variables. Reliability theory and reliability engineering also make extensive use of the exponential distribution. In physics, if you observe a gas at a fixed temperature and pressure in a uniform gravitational field, the heights of the various molecules also follow an approximate exponential distribution. In hydrology, the exponential distribution is used to analyze extreme values of such variables as monthly and annual maximum values of daily rainfall and river discharge volumes. Balakrishnan and Ahsanullah (1995) have established some recurrence relations for single and product moments of record values from exponential distribution based on record values. Ahmadi, Doostparast and Parsian (2005) have obtained the estimation and prediction based on k-record values for two-parameter exponential distribution. Ahsanullah and Kirmani (1991) have also attained some characterizations of the exponential distribution based on lower record values.

This study is going to employ a new method for finding Bayesian point and interval estimations. Consider a situation in which we only have the smallest and largest data. Consequently, the statistics based on $X_{T(1)}$ and $X_{T(n)}$ will play an important role here. This occurs in many real situations such as stock exchange. Consider a statistician who wants to make the statistical inferences about the prices of stocks and shares in a stock market, often, the middle prices are not recorded and only the highest prices and the base ones are given. Pharmacy studies could be another instance in which in order to confirm the effectiveness of drugs and poisons, the upper and lower levels of effectiveness are considered. Therefore, largest and smallest values are of great importance in decision-making.

When inference based on the largest and smallest values is required, one of the best choices is applying $R_{U.R}$ defined by $R_{U.R} = X_{T(n)} - X_{T(1)}$. In this paper, attempts have been made to draw Bayesian inferences based on this statistic, therefore, in section 2, the distribution of upper record range statistic for exponential variables has been obtained. In section 3, the Bayesian methods based on upper record range have been utilized to determine the point estimators. In section 4, credible interval with equal tails based on upper record range, and equations of Highest Posterior Density (HPD) interval have been determined. In recent years, much attention has been given to the study of the homotopy-perturbation method (HPM) to solve a wide range of problems of which mathematical models yield differential equation or system of differential equations. HPM changes difficult problems into an infinite set of problems which are now easier to solve needless to transform nonlinear terms (Momani et al. (2005), He (1999), Biazar (2007)). The applications of HPM in nonlinear problems have been demonstrated by many researchers (Ganji (2006), Dehghan et al. (2008), Chowdhury et al. (2009).). That's why it has been tried to emoloy this method to obtain HPD intervals. In section 5, the admissibility of the point estimations has been assessed.

## 2 Upper Record Range

Let $X_{T(1)}, ..., X_{T(n)}$ denote the consecutive record values observed in a sequence of independent and identically distributed random variables with probability density function f(x) and cumulative distribution function F(x). The joint pdf $x_{T(1)}, x_{T(2)}, \ldots, x_{T(n)}$ is:

$$f(x_{T(1)}, x_{T(2)}, \ldots, x_{T(n)}; \delta) = f(x_{T(n)}; \delta) \prod_{i=1}^{n-1} h(x_{T(i)}; \delta), \qquad (3)$$

where $h(x_{T(i)}; \delta) = \frac{f(x_{T(n)}; \delta)}{1 - F(x_{T(n)}; \delta)}$.

A combination of (1), (2) and (3) will lead to:

$$f(x_{T(1)}, x_{T(2)}, \ldots, x_{T(n)}; \delta) = \frac{1}{\delta^n} \exp(-\frac{x_{T(n)}}{\delta}),$$

where $x_{T(1)} < x_{T(2)} < \ldots < x_{T(n)}$.

Integrating out $x_{T(2)}, \ldots, x_{T(n-1)}$, it is gotten the joint pdf $x_{T(1)}, x_{T(n)}$ as

$$f(x_{T(1)}, x_{T(n)}) = \frac{1}{(n-2)! \delta^n} (x_{T(n)} - x_{T(1)})^{n-2} exp(-\frac{x_{T(n)}}{\delta}),$$

where

$$0 < x_{T(1)} < x_{T(n)} < \infty.$$

Using the transformations $R_{U.R} = X_{T(n)} - X_{T(1)}$, $U = X_{T(n)}$ and integrating out U, the probability density function of $R_{U.R}$ is obtained as:

$$f_{R_{U.R}}(r) = \frac{r^{n-2} exp(-\frac{r}{\delta})}{(n-2)! \delta^{n-1}}, r > 0. \tag{4}$$

This means $R_{U.R} = (X_{T(n)} - X_{T(1)})$ has been distributed as Gamma(n-1,$\delta$).

## 3 Bayesian Estimations Based on Upper Record Range Statistic

This section is going to get Bayesian estimations under three different type of loss functions (quadratic and squared error, absolute error). Our strategy depends on two conditions. First, the prior density function provided to ease the calculation of posterior distribution should be chosen. Second, since the study of the admissibility of a class of linear estimators of the form $mR_n + d$ is wanted in the next sections, the prior density function has to provide this class of estimations. The conjugate family is suitable to meet the first condition. Therefore, it is necessary to select a suitable prior distribution for $\delta$ from the conjugate family in a way that fulfils the second condition too.

with $r_n = (x_{T(n)} - x_{T(1)})$ and $L(\delta|r_n)$ as the likelihood of observing $\delta$ based on upper record range, then:

$$L(\delta|r_n) \propto \delta^{-n} \exp(\frac{-r_n}{\delta}). \tag{5}$$

According to (5), the natural conjugate prior distribution for $\delta$ can be the inverted gamma with pdf:

$$g(\delta) = \frac{b^a}{\Gamma(a) \delta^{a+1}} \exp(-\frac{b}{\delta}). \tag{6}$$

As it is shown in the next sections, this prior density function provides the expected conditions. By combining the likelihood function (5) and the prior density (6), the posterior density of $\delta$ is obtained as:

$$\pi(\delta|R_{U.R} = r) = \frac{(r+b)^{a+n-1} \exp(-\frac{r+b}{\delta})}{\Gamma(a+n-1)\delta^{a+n}}, \tag{7}$$

where

$$R_{U.R} = X_{T(n)} - X_{T(1)}, \delta > 0.$$

Note that

$$\frac{1}{\delta}|R_{U.R} = r \sim Gamma(a+n-1, \frac{1}{r+b}).$$

**Note:** From (5), logarithm of likelihood function is:

$$L(\delta; r) = (n-2)\log(r) - \log(n-2)! - (n-1)\log(\delta) - \frac{r}{\delta}, \tag{8}$$

the MLE of $\delta$ can be obtained by solving the following likelihood equation:

$$\frac{\partial L}{\partial \delta} = 0. \tag{9}$$

By solving equation (9), the MLE estimation based on upper record range for the parameter $\delta$ can be attained as:

$$\widehat{\delta}_{MBURR} = \frac{X_{T(n)} - X_{T(1)}}{n-1}. \tag{10}$$

As it is shown, this estimation is a linear of the form $mR_n + d$.

3.1 Bayesian Estimation Based On Upper Record Range Under The Quadratic Loss Function

One of the best and widely-used loss functions for estimation of scale parameters is quadratic loss function. General form of quadratic loss function is given by:

$$L = \lambda(\hat{\delta} - \delta)^2. \tag{11}$$

The quadratic loss function applied in this research is:

$$L = \frac{1}{\delta^2}(\hat{\delta} - \delta)^2. \tag{12}$$

Considering the posterior distribution based on $R_{U.R}$ (7) and Quadratic loss function (12), Bayes estimation of $\delta$ say $\widehat{\delta}_{b,1}$ (Berger 1985) will be:

$$\widehat{\delta}_{b,1} = \frac{E(\delta \frac{1}{\delta^2}|R_{U.R} = r)}{E(\frac{1}{\delta^2}|R_{U.R} = r)} = \frac{X_{T(n)} - X_{T(1)} + b}{a+n}. \tag{13}$$

According to density function (4), expectation and variance of the $\widehat{\delta}_{b,1}$ is determined as:

$$E(\widehat{\delta}_{b,1}) = \frac{(n-1)\delta + b}{a+n}, Var(\widehat{\delta}_{b,1}) = \frac{(n-1)\delta^2}{(a+n)^2}.$$

**Note:** Generally, this estimator is consistent with the parameter $\delta$. Moreover, if b=0 this estimator is asymptotically unbiased.

3.2 Bayesian Estimation Based On Upper Record Range Under The Squared Error Loss Function

Considering squared error loss function (actually in equation (11), it is assumed that $\lambda = 1$), the Bayes estimation of the parameter $\delta$ is the mean of the posterior distribution (Berger 1985), therefore:

$$\widehat{\delta}_{b,2} = E[\delta|R_{U.R}] = \frac{X_{T(n)} - X_{T(1)} + b}{a+n-2}, \quad (14)$$

and from (4)

$$E(\widehat{\delta}_{b,2}) = \frac{(n-1)\delta + b}{a+n-2}, Var(\widehat{\delta}_{b,2}) = \frac{(n-1)\delta^2}{(a+n-2)^2}.$$

**Note:** It is interesting that for a=1 and b=0, $\widehat{\delta}_{b,2}$ is transformed to an unbiased estimator of $\delta$. In fact, in this situation $\widehat{\delta}_{b,2}$ is equal to the ML estimation based on $R_{U.R}$. On the other hand $\widehat{\delta}_{b,2}$ is an asymptotically unbiased estimation for the parameter $\delta$. It is also an estimator that is in consistent with the $\delta$.

3.3 Bayesian Estimation Based On Upper Record Range Under The Absolute Error Loss Function

Another commonly-used loss function is absolute error loss function with the form of:

$$Ł(\delta, \widehat{\delta}) = |\widehat{\delta} - \delta|.$$

By considering this loss function, Bayesian estimation is the Median of the posterior distribution. According to the posterior distribution of $\delta$ (7), it is easily known that $\frac{2(R_{U.R}+b)}{\delta}$ is distributed as chi-squared with (2a+2n-2) degrees of freedom. Therefore:

$$\widehat{\delta}_{b,3} = \frac{2R_{U.R} + 2b}{\chi^2_{2n-2a-2;0.5}}.$$

Also

$$E[\widehat{\delta}_{b,3}] = \frac{2(n-1)\delta + 2b}{\chi^2_{2n-2a-2;0.5}}, Var[\widehat{\delta}_{b,3}] = \frac{4(n-1)\delta^2}{(\chi^2_{2n-2a-2;0.5})^2}.$$

## 4 Bayesian Interval estimation based on Upper Record Range Statistic

4.1 Credible Interval Based On Upper Record Range With Equal Tails

Having obtained the posterior distribution $\pi(\delta|R_{U.R})$, the problem is to see how the parameter $\delta$ is likely to lie within the interval $\mathbf{C} = [C_L, C_U]$. In the Bayesian statistics, this interval is called the credible interval. A credible interval ($\mathbf{C}$) has to meet the below condition:

$$\int_{C_L}^{C_U} \pi(\delta|R_{U.R})d\delta = 1 - \alpha.$$

Moreover, a $(1-\alpha)\%$ credible interval with equal tails must satisfy $P(\delta > b|R_{U.R} = r) = \frac{\alpha}{2}$ and $P(\delta < a|R_{U.R} = r) = \frac{\alpha}{2}$. Considering the posterior distribution (7), the below relation is obvious:

$$\frac{2}{\delta(R_{U.R}+b)}|R_{U.R} = r \sim \chi^2_{2n+2a-2}.$$

Applying the above equation, the upper and lower limits of the credible interval can be obtained from below equations:

$$P(\delta > C_U|R_{U.R}) = \frac{\alpha}{2}$$
$$P(\delta < C_L|R_{U.R}) = \frac{\alpha}{2},$$

therefore

$$P(\frac{2}{\delta(r+b)} < \frac{2}{C_U(r+b)}|R_{U.R}) = \frac{\alpha}{2} \Rightarrow \frac{2}{C_U(r+b)} = \chi^2_{2n+2a-2;\frac{\alpha}{2}},$$

and

$$P(\frac{2}{\delta(r+b)} > \frac{2}{C_L(r+b)}|R_{U.R}) = \frac{\alpha}{2} \Rightarrow \frac{2}{C_L(r+b)} = \chi^2_{2n+2a-2;1-\frac{\alpha}{2}}.$$

Hence the $(1-\alpha)\%$ credible interval with equal tails is obtained as:

$$(\frac{2}{(R_{U.R}+b)\chi^2_{2n+2a-2;1-\frac{\alpha}{2}}} < \delta < \frac{2}{(R_{U.R}+b)\chi^2_{2n+2a-2;\frac{\alpha}{2}}}).$$

Note that the length of interval is:

$$L = \frac{2}{R_{U.R}+b}(\frac{1}{\chi^2_{2n+2a-2;\frac{\alpha}{2}}} - \frac{1}{\chi^2_{2n+2a-2;1-\frac{\alpha}{2}}}).$$

4.2 Highest Posterior Density ( HPD ) Estimation Based On Upper Record Range

The Highest Posterior Density (HPD) region is defined by $\{\mathbf{C} : \pi(\delta|R_{U.R}) \geq c\}$, where $\mathbf{C} = (C_L, C_U)$ is determined from these equations:

$$\int_{c_L}^{c_U} \pi(\theta\,|R_{U.R} = r)d\theta = 1 - \alpha,$$

$$\pi(c_L\,|R_{U.R} = r) = \pi(c_U\,|R_{U.R} = r). \tag{15}$$

This interval is optimal in the sense of giving shortest lengths, for more details, see Casella and Berger (2002). By (15) and considering the posterior distribution (7) and some algebraic manipulation, the HPD equations are derived as:

$$\frac{\Gamma^*(a + n - 1, \frac{A}{c_L}, \frac{A}{c_U})}{\Gamma(a + n - 1)} = 1 - \alpha \tag{16}$$

$$(\frac{c_L}{c_U})^{a+n} = \exp\{A(\frac{1}{c_U} - \frac{1}{c_L})\},$$

where $A = b + (x_{T(n)} - x_{T(1)})$ and $\Gamma^*$ is the generalized incomplete gamma function. These equations can be solved by numerically methods.

**Note:** Length of HPD interval based on upper record range is a descending function of $\alpha$.

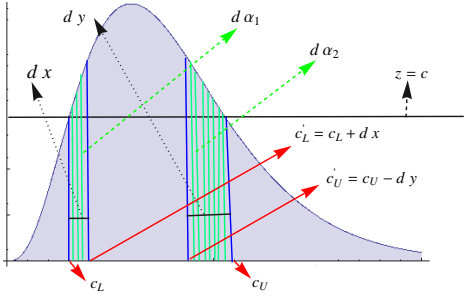

**Fig. 1** Unimodal gamma probability density function. $P(X \leq c_L) = \alpha_1$, $P(X \geq c_U) = \alpha_2$, $\alpha = \alpha_1 + \alpha_2$.

Considering figure 1

$$\mathrm{L} = C_U - C_L,$$

also

$$L^{'} = C_U - dy - (C_L - dy) = C_U - C_L - (dy + dx).$$

Therefore
$$dL = L^{'} - L = C_U - C_L - (dy + dx) - (C_U - C_L) = -(dx + dy). \quad (17)$$

On the other hand:
$$d\alpha_1 = dx \times f(c_L), d\alpha_2 = dy \times f(c_U).$$

Then
$$dx + dy = \frac{d\alpha_1}{f(c_L)} + \frac{d\alpha_2}{f(c_U)}$$

But according to (15) $f(c_L) = f(c_U)$, then
$$dx + dy = \frac{d\alpha_1 + d\alpha_2}{f(c_L)}, \quad (18)$$

by replacement (18) in (17) and noting that $d\alpha_1 + d\alpha_2 = d\alpha$
$$dL = -\frac{d\alpha}{f(c_L)}.$$

Therefore $L = c_U - c_L$ is a descending function of $\alpha$. Consequently
$$g(\alpha) = dL = B - \frac{\alpha}{f(c_L)}.$$

Also $g(\alpha)$ has to meet the below condition:
$$g(1) = 0.$$

## 5 Hight Posterior Distribution Intervals Based On Upper Record Range At The Approach Of Homotopy Perturbation Method

5.1 Basic Idea Of Homotopy Perturbation Method

Homotopy-pertuebation method (HPM) is a novel and effective method, and can solve various nonlinear equations. To illustrate the basic ideas of this method, we consider the following general nonlinear differential equation:
$$A(y) - f(r) = 0, r \in \Omega \quad (19)$$

with boundary conditions
$$B(y, \frac{\partial y}{\partial n}) = 0, r \in \Gamma, \quad (20)$$

where A is a general differential operator, B is a boundary operator, f (r) is a known analytic function , and $\Gamma$ is the boundary of the domain $\Omega$. The operator A is generally divided into two parts L and N, where L is linear while N is nonlinear. Therefore equation (19) can be written as follows:

$$L(y) + N(y) - f(r) = 0. \tag{21}$$

Now we construct a homotopy $y(r,p) : \Omega \times [0,1] \longrightarrow \Re$ of equation (19) which satisfies

$$H(y,p) = (1-p)[L(y) + L_0(y)] + p[A(y) - f(r)] = 0, p \in [0,1], r \in \Omega, \tag{22}$$

which is equivalent to

$$H(y,p) = [L(y) + L_0(y)] + pL_0 + p[N(y) - f(r)] = 0, \tag{23}$$

where $p \in [0,1]$ is an embedding parameter and $y_0$ is an initial approximation which satisfies boundary conditions. It follows from (22) and (23) that

$$H(y,0) = L(y) - L_0(y), H(y,1) = A(y) - f(r) = 0. \tag{24}$$

Thus, the changing process of p from 0 to 1 is just that of $y(r,p)$ from $y_0(r)$ to $y(r)$. In topology this is called deformation, and $L(y)L_0(y)$ and $A(y) - f(r)$ are called homotopic. Here the embedding parameter p is introduced much more naturally, unaffected by artificial factors. Due to the fact that $0 \leq p \leq 1$, so the embedding parameter can be considered as a small parameter. So it is very natural to assume that the solution of (22) and (23) can be expressed as

$$y(x) = u_0(x) + u_1(x) + u_2(x) + \cdots. \tag{25}$$

According to HPM, the approximate solution of equation (19) can be expressed as a series of the power of p, i.e.,

$$y = \lim_{p \longrightarrow 1} y = u_0 + u_1 + u_2 + \cdots. \tag{26}$$

The convergence of series (26) has been proved by He (1999).

5.2 Obtaining HPD Intervals By HPM Methods

**Theorem**: HPD intervals based on upper record range ($c_L, c_U$) with length of $c_U - c_L = g(\alpha)$ by HPM methods is obtained as follow:

$$c_L \approx \frac{A + 2ag + 2ng + \sqrt{A^2 + 8Aag + 8Ang}}{2(a+n)} \tag{27}$$

$$c_u \approx \frac{A + 2ag + 2ng + \sqrt{A^2 + 8Aag + 8Ang}}{2(a+n)} + g,$$

where $g$ can be every positive descending function from $\alpha$ providing that

$$g(0) = +\infty, g(1) = 0.$$

**proof** :

Throughout this proof, we will use the following notation:

$$g^{'} = \frac{\partial g(\alpha)}{\partial \alpha} \qquad (28)$$

HPD equations have been obtained in section (4.2) as:

$$\frac{\Gamma^{*}(a+n-1, \frac{A}{c_L}, \frac{A}{c_U})}{\Gamma(a+n-1)} = 1-\alpha, \left(\frac{c_L}{c_U}\right)^{a+n} = \exp\{A(\frac{1}{c_U} - \frac{1}{c_L})\}. \qquad (29)$$

Assume $C_L = u$ then with taking logarithm (Ln) from both sides of (29) :

$$(a+n)Lnu - (a+n)Ln(u+g) = \frac{A}{u+g} - \frac{A}{u}.$$

By taking derivation :

$$\frac{au^{'}}{u} + \frac{nu^{'}}{u} - \frac{au^{'} + ag^{'}}{u+g} - \frac{nu^{'} + ng^{'}}{u+g}$$
$$= \frac{-Au^{'} - Ag^{'}}{(u+g)^2} + \frac{Au^{'}}{u^2},$$

by simplicity

$$au^{'}u^4g + 2au^{'}u^3g^2 + au^{'}u^2g^3 + nu^{'}u^4g + 2nu^{'}u^3g^2 + nu^{'}u^2g^3$$
$$-au^5g^{'} - 2au^4g^{'}g - au^3g^{'}g^2 - nu^5g^{'} - 2nu^4g^{'}g - nu^3g^{'}g^2$$
$$= -Au^4g^{'} - Au^3g^{'}g + 2Au^{'}u^3g + 2Au^{'}u^2g^2 + Au^{'}u^2g^2 + Au^{'}ug^3,$$

or

$$(-Ag^3)u^{'} + (au^2g + 2aug^2 + ag^3 + nu^2g + 2nug^2 + ng^3 - 2Aug - 2Ag^2 - Ag^2)u^{'}u$$
$$= (-Aug^{'} - Ag^{'}g + au^2g^{'} + 2aug^{'}g + ag^{'}g^2 + nu^2g^{'} + 2nug^{'}g + ng^{'}g^2)u^2.$$

With HPM (Homotopy Perturbation Method) first construct a homotopy with embedding parameter p ($0 \leq p \leq 1$):

$$p(-Ag^3)u^{'} + p[au^2 + 2ag^2u + ag^3 + ngu^2 + 2ng^2u + ng^3 - 2Agu - 2Ag^2 - Ag]uu^{'} =$$

$$(-Ag^{'}u - Ag^{'}g + ag^{'}u^2 + 2ag^{'}gu + ag^{'}g^2 + ng^{'}u^2 + 2ng^{'}gu + ng^{'}g^2)u^2. \quad (30)$$

According to section (5.1) let u in the series form:

$$u = \sum_{+\infty}^{i=1} p^i a_i(\alpha). \qquad (31)$$

Thus:
$$u' = \sum_{i=1}^{+\infty} p^i \frac{\partial a_i(\alpha)}{\partial \alpha} = \sum_{i=1}^{+\infty} p^i a'_i, \tag{32}$$

where the series approach to the exact solutions when $p \longrightarrow 1$.
Now by substituting (31) and (32) into (30):

$$p(-Ag^3)\sum_{k=0}^{+\infty} a'_k p^k + p(a(\sum_{k=0}^{+\infty} a_k p^k)^2 + 2ag^2 \sum_{k=0}^{+\infty} a_k p^k + ag^3 \tag{33}$$

$$+ng(\sum_{k=0}^{+\infty} a_k p^k)^2 + 2ng^2(\sum_{k=0}^{+\infty} a_k p^k)^2 + ng^3 - 2Ag \sum_{k=0}^{+\infty} a_k p^k -$$

$$2Ag^2 - Ag^2)\sum_{k=0}^{+\infty} a'_k p^k \sum_{i=0}^{+\infty} a_i p^i$$

$$= (-Ag'\sum_{k=0}^{+\infty} a_k p^k - Ag'g - ag'(\sum_{k=0}^{+\infty} a_k p^k)^2 + 2ag'g(\sum_{k=0}^{+\infty} a_k p^k) +$$

$$+ag'g^2 + ng'(\sum_{k=0}^{+\infty} a_k p^k)^2 + 2ng'g \sum_{k=0}^{+\infty} a_k p^k + ng'g^2)(\sum_{k=0}^{+\infty} a_k p^k)^2.$$

Collecting the terms of both sides of (33) with the same powers of $p$ ($p^0, p^1, p^2, \cdots$) and equating each coefficient of $p$ from both sides of (33) results in: (by this algorithm $a_0, a_1, a_2, \cdots$ can be determined)

$$p^0 : (-Ag'a_0 - Ag'g + ag'a_0^2 + 2ag'ga_0 \tag{34}$$
$$+ag'g^2 + ng'a_0^2 + 2ng'ga_0 + ng'g^2)a_0^2 = 0.$$

$$p^1 : (-Ag^3 a'_0 + aa_0^2 + 2ag^2 a_0 + ag^3 + nga_0^2$$
$$+2ng^2 a_0^2 + ng^3 - 2Aga_0 - 2Ag^2 - Ag^2)a'_0 a_0 =$$
$$-3Ag'a_0^2 a_1 - 2Ag'ga_0 a_1 + 4ag'a_0^3 a_1 + 6ag'ga_0^2 a_1$$
$$+2ag'g^2 a_0 a_1 + 4ng'a_0^3 a_1 + 6ng'ga_0^2 a_1 + 2ng'g^2 a_0 a_1.$$

Now (34) is a second order equation ($Ma_0^2 + Na_0 + Q = 0$) and by solving it $a_0(\alpha)$ be determined:

$$(a+n)a_0^2 + (A + 2ag + 2ng)a_0 + (ag^2 - Ag + ng^2) = 0$$
$$\triangle = N^2 - 4MQ > 0,$$

thus, $a_0$ can be calculated as follows :

$$a_0 = \frac{A + 2ag + 2ng + \sqrt{A^2 + 8Aag + 8Ang}}{2(a+n)}.$$

**Remark** :
Note $a_0$ satisfies (34) too, but when $a_0=0$ then $a_1, a_2, \cdots =0$. Therefore, $a_0 = 0$ cannot be considered.

Now $c_L = u = \sum_{i=0}^{+\infty} p^i a_i(\alpha)$ and then:

$$c_L = u \approx a_0(\alpha) = \frac{A + 2ag + 2ng + \sqrt{A^2 + 8Aag + 8Ang}}{2(a+n)},$$

and by $c_u = c_L + g(\alpha)$:

$$c_U = u \approx \frac{A + 2ag + 2ng + \sqrt{A^2 + 8Aag + 8Ang}}{2(a+n)} + g(\alpha).$$

## 6 Admissibility and point estimations

Attention should be paid that the obtained estimators in the section 3 are all a special form of the linear model $mR_n + d$ (or $m(X_{T(n)} - X_{T(1)}) + d$). This paves the way for further studies. In this part, it is tried to determine the range of m and d by which the linear form will have desirability. One of the criteria which has been considered in statistical inference and in decision theory, for this aim, is the admissibility concept. To find the fundamental arguments relating to admissibility, one may refer to Blyth (1951).

**Theorem**: Estimators with the form of $m(X_{T(n)} - X_{T(1)}) + d$ under quadratic loss function ($L = \frac{1}{\delta^2}(\hat{\delta} - \delta)^2$) are admissible estimators if:
1- $m\epsilon[0, 1/n)$ and $d > 0$
2- $m = 1/n$ and $d > 0$

**Proof** :
First, due to quadratic loss function (12) and considering probability density function of upper record range (4) and prior density function with the form of (6), risk function and Bayesian risk will be obtained in the following way:

$$R(mR_n + d, \delta) = [(m(n-1) - 1)^2 + m^2(n-1)] + \frac{[2d(m(n-1)-1)]}{\delta} + \frac{d^2}{\delta^2}$$

$$\begin{aligned}r(mR_n + d, \delta) = &[(m(n-1) - 1)^2 + m^2(n-1)] + \\ &\frac{2da[m(n-1) - 1]}{b} + \frac{d^2 a(a+1)}{b^2}\end{aligned} \quad (35)$$

**1:** Due to Bayesian estimator under quadratic loss function ($L = \frac{1}{\delta^2}(\hat{\delta} - \delta)^2$), and noticing that Bayesian estimators always have the quality of admissibility, and also paying attention to:

$$\hat{\delta} = \frac{R_n + b}{a + n} = \frac{R_n}{a + n} + \frac{b}{a + n} = mR_n + d$$

$$a \to +\infty \Rightarrow \frac{1}{a + n} \to 0 \Rightarrow m \to 0$$

$$a \to 0 \Rightarrow \frac{1}{a + n} \to \frac{1}{n} \Rightarrow m \to \frac{1}{n} \tag{36}$$

$$a \to 0 \Rightarrow \frac{b}{a + n} \to \frac{b}{n} \Rightarrow d \to \frac{b}{n} \tag{37}$$

then linear estimator of $mR_n + d$ has the quality of admissibility for $m\epsilon(0, 1/n)$ and $b > 0$.

Also by considering (35) and assuming that $m = 0$, we have $R(mR_n + d, \delta) = \frac{1}{\delta^2}(\delta - d)^2$. Therefore, if $\delta = d$, it is seen that $R = 0$. It means $mR_n + d$ is admissible for $m = 0$ and $d > 0$.

**2:** By designating a as $1/k$ in prior density function, then: (6)

$$g(\delta) = \frac{b^{\frac{1}{k}}}{\Gamma(\frac{1}{k})\delta^{\frac{1}{k}+1}} \exp(-\frac{b}{\delta}). \tag{38}$$

Also if

$$m = \frac{k}{1 + kn}, d = \frac{kb}{1 + kn},$$

Regarding the first section of the theorem, $\hat{\delta} = \frac{kR_n}{1+kn} + \frac{kb}{1+kn}$ is an admissible estimator

$$0 \leq m = \frac{k}{1 + kn} < \frac{1}{n}, 0 < d = \frac{kb}{1 + kn}.$$

Considering the mentioned factors and prior density function (38), Bayesian risk function will be obtained in the following way:

$$r_1 = [[((\frac{k}{1 + kn})(n - 1)) - 1]^2 + (\frac{k}{1 + kn})^2(n - 1)] +$$

$$\frac{2\frac{kb}{1+kn}[(\frac{k}{1+kn})(n - 1) - 1]}{kb} + \frac{k + 1}{k^2}\frac{(\frac{kb}{1+kn})^2}{b^2}. \tag{39}$$

On the other hand, by assuming the linear estimator in the form of $\hat{\delta} = (\frac{R_n}{n}) + (\frac{1}{n})$, its related Bayesian risk function under prior density function (38) would be gained as follows:

$$r_2 = \frac{1}{n} - \frac{2}{bn^2k} + \frac{k + 1}{n^2k^2b^2}. \tag{40}$$

Considering (39) and (40)

$$\lim(r_1 - r_2)_{k \to +\infty} = 0$$

Therefore, $\hat{\delta} = (\frac{R_n}{n}) + (\frac{1}{n})$ has the admissibility quality because it will be impossible to have one other $(\hat{\delta}')$ as given in $Risk(\hat{\delta}') < Risk(\hat{\delta})$, otherwise the following contradiction will occur:

$$Risk(\hat{\delta}') < Risk(\hat{\delta}) \Rightarrow r(\hat{\delta}') < r_1(\hat{\delta}) = r_2(\hat{\delta})$$

**Theorem**: Estimators with the form of $m(X_{T_n} - X_{T_1}) + d$ under squared error loss function $L = (\hat{\delta} - \delta)^2$ are admissible estimators if:
1- $m \epsilon [0, 1/n), d > 0$
2- $m = \frac{1}{n}, d > 0$.

**Proof** :
The procedures of proof is exactly similar to the one previously mentioned.

## 7 Numerical Results

7.1 Simulation without repetition

**Example1** :
To illustrate the developed estimation techniques, consider the following simulated data from the exponential distribution:

0.06274109, 4.38197283, 5.64659541, 0.08382565, 5.27747401, 2.69666048, 0.98792501, 2.36520919, 0.04765528, 0.63918881, 0.07107701, 2.19439004, 2.71178500, 1.45946486, 5.31182137, 0.42911833, 2.74980209, 0.41108542, 2.21423065, 1.31309101, 0.29502675, 1.50707359, 7.26620864, 2.47032883, 2.79500172, 1.14469466, 3.20462205, 4.10787212, 2.97814895, 2.42587180, 1.85331396, 0.70619791, 2.60601466, 1.28472926, 0.29126746, 0.07298126, 0.24644642, 1.90989237, 2.40637729, 2.17449704, 1.02288571, 1.54665282, 2.95083160, 0.95526777, 0.04135414, 1.01268457, 1.07257669, 0.75808989, 3.33255820, 0.71060492, 1.18752218, 9.41371352, 9.51953091.

This data has been obtained using the transformation $X_i = (-\delta log(1 - U_i))$ where $U_i$ is a uniformly distributed random variable. Upper record values from this sample have been observed as: 0.06274109, 4.38197283, 5.64659541, 7.26620864, 9.41371352, 9.51953091.
Consequently, the upper record range values are obtained as: 4.319232, 5.583854, 7.203468, 9.350972, 9.456790. For prior distribution with parameters ($a = 3, b = 5$), Bayesian and ML estimations, related MSE's based on upper record range ((10), (13) and (14)) have been also computed for n=2, 3, 4, 5, 6 (Table1). Moreover, through the table 1 and figure 2, estimations based on upper record range have been compared with both ML estimation based on records and ML estimation based on sample by MSE's. As it is shown in this simulation, generally, Bayesian estimations based on upper record range have less MSEs than ML estimations based on the records and the sample.

Table 1. Estimations and MSE's

|   | Estimator | Estimated value | MSE |
|---|---|---|---|
| n=2 | $\hat{\theta}_{MLE}$ | 2.975741 | 2.231147 |
| 3 |  | 3.866026 | 1.841971 |
| 4 |  | 2.543784 | 0.7983858 |
| 5 |  | 3.090522 | 0.6315874 |
| 6 |  | 3.024878 | 0.3991685 |
| n=2 | $\hat{\theta}_{MLE-based-on-records}$ | 2.19098642 | 2.400210736 |
| 3 |  | 1.88219847 | 1.180890361 |
| 4 |  | 1.81655216 | 0.824965437 |
| 5 |  | 1.88274270 | 0.708944018 |
| 6 |  | 1.58658848 | 0.419543836 |
| n=2 | $\hat{\theta}_{MLE-}$ | 4.319232 | 18.655763 |
| 3 | based-on-upper-record-range | 2.791927 | 3.897429 |
| 4 |  | 1.81655216 | 0.824965437 |
| 5 |  | 2.337743 | 1.366261 |
| 6 |  | 1.891358 | 0.715447 |
| n=2 | $\hat{\theta}_{b,1}$ | 1.863846 | 0.3157106 |
| 3 |  | 1.763976 | 0.2154080 |
| 4 |  | 1.743353 | 0.1862123 |
| 5 |  | 1.793872 | 0.1860660 |
| 6 |  | 1.606310 | 0.1242328 |
| n=2 | $\hat{\theta}_{b,2}$ | 3.106411 | 0.1280480 |
| 3 |  | 2.645964 | 0.1257608 |
| 4 |  | 2.440694 | 0.1200945 |
| 5 |  | 2.391829 | 0.1113384 |
| 6 |  | 2.065256 | 0.1056581 |

**Example2** : As another example, below data has been simulated from exponential distribution with ($\delta = 1$), for $a = 3$ and $b = 4$, the shortest interval and credible interval with equal tails based on upper record range have been obtained for $\delta$ (Table 2, 3).

0.067773, 0.056655, 0.032254, 2.081551, 0.125478, 2.002154, 1.9874521, 1.254875,
0.236587, 1.876541, 0.231456, 2.274237, 0.336521, 1.985436, 2.001245, 3.373468,
2.125987, 1.236541, 0.236541, 1.789654, 3.021543, 2.365987, 1.002154, 0.357951,
2.147963, 3.123623, 2.543659, 1.598723, 0.001357, 1.986124, 1.963254, 3.847746,
2.356547, 1.235463, 1.4723568, 1.983217, 3.002541, 4.243143.

Table 2. Shortest interval estimation and length .

| Number of record | Confidence level | Lower limit $C_L$ | Upper limit $C_U$ | Length |
|---|---|---|---|---|
| n=2 | 90% | 0.2374295 | 20.321103 | 20.083673 |
| 3 | | 0.2437247 | 4.178513 | 3.934789 |
| 4 | | 0.3326670 | 3.039720 | 2.707053 |
| 5 | | 0.3430025 | 2.206510 | 1.863508 |
| 6 | | 0.3429490 | 1.754947 | 1.411998 |
| n=2 | 95% | 0.1990373 | 41.712780 | 41.513743 |
| 3 | | 0.2095089 | 6.231189 | 6.021680 |
| 4 | | 0.2907256 | 5.803538 | 5.512812 |
| 5 | | 0.3034332 | 2.797183 | 2.493750 |
| 6 | | 0.3062671 | 2.149206 | 1.842939 |
| n=2 | 99% | 0.1467439 | 212.714428 | 212.567684 |
| 3 | | 0.1607021 | 14.863348 | 14.702646 |
| 4 | | 0.2288248 | 11.946856 | 11.718031 |
| 5 | | 0.2434710 | 4.610567 | 4.367096 |
| 6 | | 0.2495104 | 3.283298 | 3.033788 |

Table 3. Credible interval estimation with equal tails and length.

| Number of record | Confidence level | Lower limit $C_L$ | Upper limit $C_U$ | Length |
|---|---|---|---|---|
| n=2 | 90% | 0.07174895 | 1.55242910 | 1.48068016 |
| 3 | | 0.05359375 | 1.30413666 | 1.25054291 |
| 4 | | 0.03116289 | 1.14478602 | 1.11362313 |
| 5 | | 0.02389973 | 1.09615636 | 1.07225663 |
| 6 | | 0.01929951 | 1.06956777 | 1.05026826 |
| n=2 | 95% | 0.06252402 | 1.73013885 | 1.66761482 |
| 3 | | 0.04739758 | 1.38128334 | 1.33388576 |
| 4 | | 0.02785230 | 1.17570219 | 1.14784989 |
| 5 | | 0.02153339 | 1.11410992 | 1.09257654 |
| 6 | | 0.01750089 | 1.08120917 | 1.06370827 |
| n=2 | 99% | 0.41941285 | 2.3369799 | 1.91756703 |
| 3 | | 0.24271949 | 1.6181843 | 1.37546478 |
| 4 | | 0.11919901 | 1.2646281 | 1.14542911 |
| 5 | | 0.08088953 | 1.1634824 | 1.08259283 |
| 6 | | 0.05946667 | 1.1121817 | 1.05271508 |

## 7.2 Simulation with repetition

In the following example it has been tried to show simulation with repetition (number of repetition is 100) for different parameters a and b. Same as before subsection data has been obtained using the transformation $X_i = (-\delta log(1 - U_i))$ where u is a uniformly distributed random variable.

Table 4. Estimations and MSE's (N.R=100 and a=8, b=2, $\delta = 2$)

| | Estimator | Average of estimated value | MSE |
|---|---|---|---|
| n=4 | $\hat{\theta}_{MLE-Based-On-Upper-Records}$ | 0.4588271 | 0.089599176 |
| 5 | | 0.5000819 | 0.075141503 |
| 6 | | 0.4768701 | 0.054178733 |
| 7 | | 0.4967541 | 0.047924891 |
| n=4 | $\hat{\theta}_{b,1}$ | 0.2813734 | 0.005742126 |
| 5 | | 0.3077175 | 0.008148448 |
| 6 | | 0.3131679 | 0.008353949 |
| 7 | | 0.332035 | 0.009285333 |
| n=4 | $\hat{\theta}_{b,2}$ | 0.3376481 | 0.007469059 |
| 5 | | 0.3636661 | 0.009756746 |
| 6 | | 0.3653625 | 0.009625019 |
| 7 | | 0.3831173 | 0.009869343 |

Table 5. Estimations and MSE's (N.R=100 and a=6, b=1, $\delta = 3$)

| | Estimator | Average of estimated value | MSE |
|---|---|---|---|
| n=4 | $\hat{\theta}_{MLE-Based-On-Upper-Records}$ | 0.3172188 | 0.043272931 |
| 5 | | 0.3331795 | 0.034504358 |
| 6 | | 0.326504 | 0.024647407 |
| 7 | | 0.3467425 | 0.023494405 |
| n=4 | $\hat{\theta}_{b,1}$ | 0.1951656 | 0.003849557 |
| 5 | | 0.2120653 | 0.004990664 |
| 6 | | 0.2193766 | 0.004746278 |
| 7 | | 0.2369581 | 0.005640035 |
| n=4 | $\hat{\theta}_{b,2}$ | 0.243957 | 0.005340475 |
| 5 | | 0.2591909 | 0.006308964 |
| 6 | | 0.263252 | 0.005712699 |
| 7 | | 0.2800414 | 0.006093913 |

## 8 Conclusion

In this paper, Bayesian point estimations based on upper record range have been obtained for three different loss functions. It has been revealed that they hold interesting properties such as consistency, asymptotically unbiaseness, and admissibility. Moreover, it was observed that all of them take the form of $mR_n + d$. Also, general admissibility conditions for this form has been obtained by means of two theorems. This appears to be appealing and verifiable because by applying only two values from the upper records, it was possible to acquire other estimations which have minimum risk functions (admissibility). Additionally, Bayesian interval estimations based on upper record range have been obtained. It was shown that the length of HPD interval based on upper record is an increasing function of confidence level $(1-\alpha)$. The above function was also obtained. By means of this function and a new numerical method called HPD, limits were obtained in general form. Finally, through simulation, it was found that the MSE's of Bayesian estimations based on upper record range are less than both MSE of MLE based on records and MSE of MLE based on sample. Also, it was revealed th at credible interval based on upper record range with equal tails has a shorter length from the shortest interval estimation.

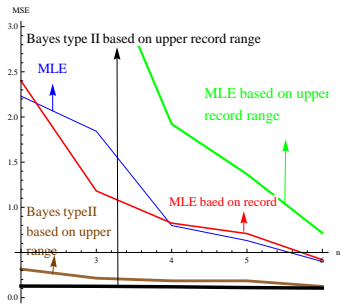

**Fig. 2** MSE's of the estimators $\hat{\theta}_{MLE}$ , $\hat{\theta}_{b,1}$ and $\hat{\theta}_{b,2}$ .

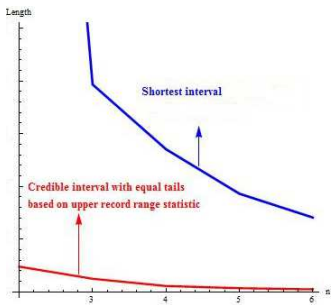

**Fig. 3** Comparison of lengths for 90% confidence.

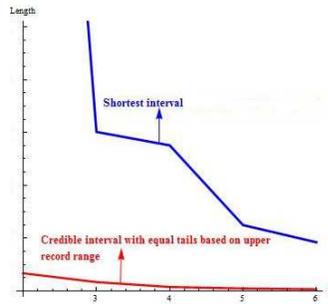

**Fig. 4** Comparison of lengths for 95% confidence.

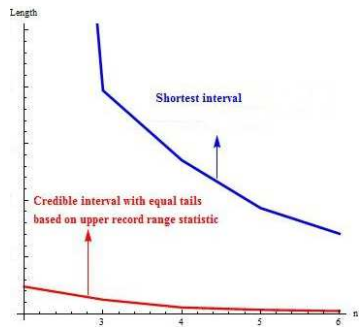

**Fig. 5** Comparison of lengths for 99% confidence.

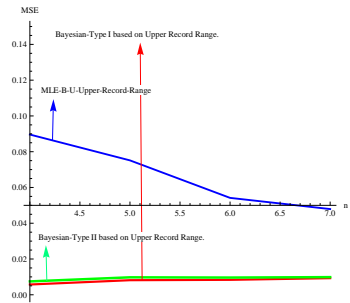

**Fig. 6** MSE's of the estimators ($\hat{\theta}_{MLEBBURR}$, $\hat{\theta}_{b,1}$ and $\hat{\theta}_{b,2}$ N.R=100 and a=8, b=2, $\delta = 2$.

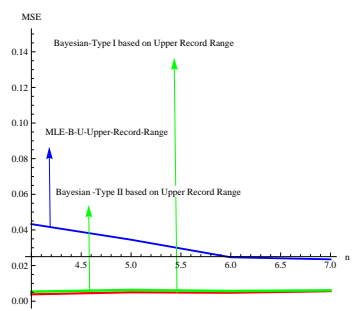

**Fig. 7** MSE's of the estimators ($\hat{\theta}_{MLEBBURR}$ , $\hat{\theta}_{b,1}$ and $\hat{\theta}_{b,2}$  N.R=100 and a=6, b=1, $\delta = 3$ .